
\input amstex
\documentstyle{amsppt}

\magnification=1200
\parskip 5pt
\pagewidth{5.4in} \pageheight{7.1in}

\baselineskip=14pt \expandafter\redefine\csname
logo\string@\endcsname{} \NoBlackBoxes \NoRunningHeads

\redefine\qed{$\ \blacksquare$}

\def\r1{\sqrt{-1}}

\def\t{\frak{t}}

\def\r1{\sqrt{-1}}

\def\bR{\mathord{\Bbb R}}

\def\bC{\mathord{\Bbb C}}

\topmatter
\title Polynomial representatives of Schubert classes in $QH^*(G/B)$
\endtitle
\author Augustin-Liviu Mare
\endauthor
\abstract We show how the quantum Chevalley  formula for $G/B$, as
stated by Peterson and   proved rigorously by Fulton and Woodward,
combined with ideas of Fomin, S. Gelfand and Postnikov, leads to a
formula which describes polynomial representatives of the Schubert
cohomology classes in the canonical presentation of $QH^*(G/B)$ in
terms of generators and relations. We generalize in this way
results of \cite{FGP}.
\endabstract
\endtopmatter
\document

\head \S 1 Introduction
\endhead
A  theorem of Borel \cite{B} describes  the
cohomology\footnote{The coefficient ring for  cohomology will
always be $\bR$.} ring of the generalized complex flag manifold
$G/B$ as the ``co-invariant algebra" of the Weyl group of $G$,
which is essentially a quotient of a certain polynomial ring. The
Schubert cohomology classes (i.e. Poincar\'e duals of Schubert
varieties) are a basis of $H^*(G/B)$. In order to determine the
structure constants of the cup-multiplication on $H^*(G/B)$ with
respect to this basis, we need to describe the Schubert cohomology
classes in Borel's presentation. According to Bernstein, I. M.
Gelfand and S. I. Gelfand [BGG], we obtain polynomial
representatives of Schubert classes in Borel's ring by starting
with a representative of the top cohomology and then applying
successively {\it divided difference operators} associated to the
simple roots of $G$. More details concerning the
Bernstein-Gelfand-Gelfand construction can be found in section 2
of our paper.

When dealing with the (small) quantum cohomology ring $QH^*(G/B)$
we face a similar situation. There exists a canonical presentation
of that ring,  again as a  quotient of a polynomial algebra, where
the variables are the same as in the classical case, plus the
``quantum variables" $q_1,\ldots , q_l$. As about the ideal of
relations, it is generated by the ``quantum deformations" of the
relations from Borel's presentation of $H^*(G/B)$ (for more
details, see section 3). The Schubert classes are a basis of
$QH^*(G/B)$ as a $\bR[q_1,\ldots,q_l]$-module.  A natural aim (see
the next paragraph) is to describe them in the previous
presentation of $QH^*(G/B)$. Our main result gives a method for
obtaining such polynomial representatives.  It can be described
briefly as follows: we start with an {\it arbitrary} polynomial
representing the Schubert class $\sigma_w$ in Borel's description
(e.g. by using the B-G-G construction); this is transformed into a
polynomial representing $\sigma_w$ in the canonical description of
$QH^*(G/B)$  after successive applications of divided difference
operators, multiplications by $q_j$'s and integer numbers and
additions. The precise formula is stated in Theorem 3.6 (also see
Lemma 3.4 and relation (5) in order to understand the notations).
For  $G=SL(n,\bC)$, the same result was proved by Fomin, S.
Gelfand and Postnikov \cite{FGP}. The main ingredient of our proof
is a result of D. Peterson \cite{P} (we call it the ``quantum
Chevalley formula", since Chevalley obtained a similar result for
the cup-multiplication on $H^*(G/B)$) which describes the quantum
multiplication by degree 2 Schubert classes.

Finally, a few words should be said about the importance of our
result. The standard presentation of $QH^*(G/B)$ mentioned above
is explicitly determined by Kim \cite{K} (see also \cite{M}). Our
description  could be relevant for finding the  structure
constants of the quantum multiplication with respect to the basis
consisting of Schubert classes, which would lead immediately to
the Gromov-Witten invariants of $G/B$. The efficiency of this
strategy depends very much on the input: we dispose of the  {\it
choice} of  polynomial representatives of Schubert classes in
Borel's ring and this has to be made judiciously (see again
\cite{FGP}, as well as Billey and Haiman \cite{BH} and Fomin and
Kirillov \cite{FK}).

\head \S 2 The Bernstein-Gelfand-Gelfand construction
\endhead
The main object of study of this paper is the generalized
 complex flag manifold $G/B$, where $G$ is a connected, simply
connected,  semisimple, complex Lie group and $B\subset G$ a Borel
subgroup. Let  $\t$ be the Lie algebra of a maximal torus of a
compact real form of $G$ and  $\Phi \subset \t^*$ the
corresponding set of roots. The negative of the Killing form
restricted to $\t$ gives an inner product $\langle \ , \ \rangle$.
To any root $\alpha$ corresponds the coroot
$$\alpha^{\vee}:= \frac{2\alpha}{\langle \alpha,
\alpha\rangle}$$ which {\it is an element of $\t$}, by using the
identification of $\t$ and $\t^*$ induced by $\langle \ , \
\rangle$. If $\{\alpha_1, \ldots ,\alpha_l\}$ is a system of
simple roots then $\{\alpha_1^{\vee},\ldots, \alpha_l^{\vee}\}$ is
a  system of simple coroots.  Consider  $\{\lambda_1 ,\ldots ,
\lambda_l\} \subset \t^*$ the corresponding system of fundamental
weights, which are defined by
$\lambda_i(\alpha_j^{\vee})=\delta_{ij}$.  To any positive root
$\alpha$ we assign the reflection $s_{\alpha}$ of $(\t, \langle \
, \ \rangle)$ about the hyperplane $\ker \alpha$. The Weyl group
$W$ is generated by all reflections $s_{\alpha}$, $\alpha\in
\Phi^+$: it is actually generated by a smaller set, namely by the
simple reflections $s_1=s_{\alpha_1}, \ldots , s_l=s_{\alpha_l}$.
To any $w\in W$ corresponds  a length, $l(w)$, which is the
smallest number of factors in a decomposition of $w$ as a product
of simple reflections.

There are two different ways to describe $H^*(G/B)$: On the one
hand, we can take $B^-\subset G$ the Borel subgroup opposite to
$B$ and assign to each $w\in W$ the Schubert variety
$\bar{C}_w=\overline{B^-.w}$, which has real codimension $2l(w)$;
its Poincar\'e dual $\sigma_w$ is  an element of $H^{2l(w)}(G/B)$;
the set $\sigma_w$, $w\in W$ is a basis of $H^*(G/B)$.  On the
other hand, let us  consider the symmetric algebra  $S(\t^*)$,
which consists of polynomial functions on $\t$. A theorem of Borel
 says that the ring homomorphism $S(\t^*)\to H^*(G/B)$
induced by $\lambda_i\mapsto \sigma_{s_i},$ $1\leq i \leq l$, is
surjective; moreover it induces the  ring isomorphism
$$ H^*(G/B) \simeq \bR[\lambda]/I_W, \tag 1$$ where $I_W$ is the
ideal of $S(\t^*)=\bR[\lambda_1,\ldots,\lambda_l]=\bR[\lambda]$
generated by the $W$-invariant polynomials of strictly positive
degree.

One is looking for a Giambelli type formula, which connects these
two descriptions by assigning to each Schubert cycle $\sigma_w$ a
polynomial representative in the quotient ring $\bR[\lambda]/I_W$.
We are going to sketch  the construction of such polynomials, as
performed by Bernstein, I. M. Gelfand and S. I. Gelfand in [BGG].
It relies on the following facts: \item {$\bullet$}  $H^*(G/B)$
and $\bR[\lambda]/I_W$ are generated as rings by $\sigma_{s_i}$,
respectively $\lambda_i$, $1\leq i \leq l$, \item {$\bullet$}
 we have a formula of Chevalley which gives the matrix of the
cup multiplication by $\sigma_{s_i}$ on $H^*(G/B)$ with respect to
the basis $\{\sigma_w:w\in W\}$, \item {$\bullet$}  there is
another, ``very similar", formula, which involves the divided
difference operators $\Delta_w$, $w\in W$ (see below) on the
polynomial ring $\bR[\lambda]$.

The following result was  proved by Chevalley  [Ch] (see also
Fulton and Woodward [FW]). \proclaim{Lemma 2.1}{\rm (Chevalley's
formula)}. For any $1\leq i \leq l$ and any $w\in W$ we have
$$\sigma_{s_i}\sigma_w=\sum_{\alpha\in \Phi^+, l(ws_{\alpha})=l(w)+1}
\lambda_i(\alpha^{\vee})\sigma_{ws_{\alpha}}.$$
\endproclaim

To each positive root $\alpha$ we assign the {\it divided
difference operator} $\Delta_{\alpha}$ on the ring $\bR[\lambda]$
(the latter being just the symmetric ring $S(\t^*)$,  it admits a
natural action of the Weyl group $W$):
$$\Delta_{\alpha}(f)=\frac{f-s_{\alpha}f}{\alpha}$$
If $w$ is an arbitrary element of  $W$, take $w =s_{i_1} \ldots
s_{i_k}$  a reduced expression and then set $$\Delta _w= \Delta
_{\alpha _{i_1}} \circ \cdots \circ \Delta_{\alpha _{i_k}}.$$ One
can show (see for instance [Hi]) that the definition does not
depend on the choice of the reduced expression. The operators
obtained in this way have the following property:
$$\Delta _w
\circ \Delta _{w'} =\cases
                                        \Delta_{ww'},&\text{ if \ }
                                        l(ww')=l(w)+l(w')\\
                                        0,&\text{ otherwise}.
                                        \endcases
\tag  2 $$

The importance of those operators for our present context is
revealed by the similarity of the following formula with Lemma
2.1: \proclaim{Lemma 2.2} {\rm (Hiller [Hi])} If $\lambda_i^*$
denotes the operator of  multiplication by $\lambda_i$ on
$\bR[\lambda]$, then for any $w\in W$ we have
$$\Delta_w\lambda_i^*-w\lambda_i^* w^{-1}\Delta_w=
\sum_{\beta\in \Phi^+,l(ws_{\beta})=l(w)-1}\lambda_i(\beta^{\vee})
\Delta_{ws_\beta}.$$
\endproclaim

Let $w_0$ be the longest element of $W$. The polynomial
$$c_{w_0}:=\frac{1}{|W|}\prod_{\alpha \in \Phi^+}\alpha$$ is
homogeneous, of degree $l(w_0)$   and
 has the  property that $\Delta_{w_0}c_{w_0}=1$. But $l(w_0)$ is
 at the same time the complex dimension of $G/B$, and it can be easily
 shown  that the class of $c_{w_0}$ in $\bR[\lambda]/I_W$ generates the top cohomology
  of $G/B$.
 To any $w\in W$ we assign $c_w:= \Delta_{w^{-1}w_0}c_{w_0}$
which is a homogeneous polynomial  of degree
$l(w)$ satisfying
$$\Delta_v c_w=  \cases
   c_{wv^{-1}}, & \text{if \ } l(wv^{-1})=l(w)-l(v) \\
    0, & \text{otherwise}
  \endcases
$$
for any $v\in W$ (see (2)). In particular, if $l(v)=l(w)$, then
$\Delta_v(c_w)=\delta_{vw}$. Since  $\Delta_w$ leaves $I_W$
invariant, it induces an operator on $\bR[\lambda]/I_W$ which also
satisfies $\Delta_v([c_w])=\delta_{vw}$, provided that
$l(v)=l(w)$. Because $\dim \bR[\lambda]/I_W=|W|$, it follows that
the classes $[c_w]$, $w\in W$, are a basis of $\bR[\lambda]/I_W$.
We can easily determine any of the coefficients $a_v$ from
$$\lambda_i[c_w]=\sum_{l(v)=l(w)+1}a_v[c_v],$$
by applying $\Delta_v$ on both sides and using Lemma 2.2. It
follows that
$$
\lambda_i [c_w]=\sum_{\alpha\in \Phi^+, l(ws_{\alpha})=l(w)+1}
\lambda_i(\alpha^{\vee})[c_{ws_{\alpha}}]. \tag 3 $$ From
$\Delta_{s_i}(\lambda_j)=\delta_{ij}$, $1\leq i,j\leq l$, we
deduce that $c_{s_i}=\lambda_i$. We just have to compare (3) with Lemma
2.1 to conclude: \proclaim{Theorem 2.3}{\rm  (Bernstein, I. M.
Gelfand and  S. I. Gelfand \cite{BGG})} Let $[c_{w_0}]$ be the
image of $\sigma_{w_0}$ by the identification
$H^*(G/B)=\bR[\lambda]/I_W$ indicated above. Then the map
$\sigma_w\mapsto [c_w]:=\Delta_{w^{-1}w_0}[c_{w_0}]$ is a ring
isomorphism.
\endproclaim
The polynomial $c_w=\Delta_{w^{-1}w_0}c_{w_0}$ being a
representative of the Schubert cycle $\sigma_w$ in
$\bR[\lambda]/I_W$, is a solution of the classical (i.e.
non-quantum) Giambelli problem for $G/B$.

\head \S 3 Quantization map
\endhead

Additively, the quantum cohomology  $QH^*(G/B)$ of $G/B$ is just
$H^*(G/B)\otimes \bR[q_1,\ldots, q_l]$, where $l$ is the rank of
$G$ and $q_1, \ldots, q_l$ are some  variables. The multiplication
$\circ$ is uniquely determined by $\bR[q]$-linearity and the
general formula
$$ \sigma_u \circ \sigma_v =\sum_{d=(d_1,\ldots ,d_l)\geq 0}
q^d\sum_{w\in W} \langle \sigma_u|\sigma_v|\sigma_{w_0w}\rangle_d
\sigma_w,$$ $u,v\in W$, where $q^d$ denotes $q_1^{d_1}\ldots
q_l^{d_l}$. The coefficient
 $\langle \sigma_u|\sigma_v|\sigma_{w_0w}\rangle_d$ is the
Gromov-Witten invariant, which counts the number of holomorphic
curves $\varphi :\bC P^1 \to G/B$ such that $\varphi_*([\bC P^1])=
d$ in $H_2(G/B)$ and $\varphi(0)$, $\varphi(1)$ and
$\varphi(\infty)$ are in general translates of the Schubert
varieties dual to    $\sigma_u$, $\sigma_v$, respectively
$\sigma_{w_0w}$. It turns out that  this number can be nonzero and
finite only if $l(u)+l(v)=l(w)+2\sum_{i=1}^l d_i$; if it is
infinity, we set $ \langle
\sigma_u|\sigma_v|\sigma_{w_0w}\rangle_d=0$. The ring
$(QH^*(G/B),\circ)$ is  commutative and associative (for more
details about quantum cohomology we refer the reader to Fulton and
Pandharipande \cite{FP}).

One can show that the quantum cohomology ring of $G/B$ is
generated by $H^2(G/B)\otimes \bR[q_1,\ldots, q_l]$, i.e. by $q_1,
\ldots , q_l,\lambda_1,\ldots,\lambda_l$. To determine the ideal
of relations, we only have to take any of the fundamental
$W$-invariant polynomials $u_i$, $1\leq i \leq l$
--- as generators of the ideal $I_W$ of relations in $H^*(G/B)$
--- and find its ``quantum deformation" $R_i$. The latter is a
polynomial in $\bR[q,\lambda]$, uniquely determined by:

(a) the relation $R_i(q_1,\ldots ,q_l,\sigma_{s_1}\circ,\ldots,
\sigma_{s_l}\circ)=0$ holds in $QH^*(G/B)$,

(b) the component of $R_i$ free of $q$ is $u_i$.

\noindent If $I_W^q$ denotes the ideal of $\bR[q,\lambda]$
generated by $R_1, \ldots, R_l$, then we have the ring isomorphism
$$
QH^*(G/B)\simeq \bR[q,\lambda]/I_W^q. \tag 4 $$

 The challenge
is now to solve the ``quantum Giambelli problem": via the
isomorphism (4),
 find a polynomial representative in $\bR[q,\lambda]/I_W^q$
for each Schubert class $\sigma_w$, $w\in W$. We can actually use
Theorem 2.3 in order to rephrase the problem as follows: Describe
(the image of $[c_w]$ via) the map
 $$
\bR[q,\lambda]/(I_W\otimes\bR[q])=\bR[\lambda]/I_W\otimes \bR[q]
@>{\cong}>>H^*(G/B)\otimes \bR[q] =QH^*(G/B)@>{\cong}>>
\bR[q,\lambda]/I_W^q.
$$
 Note that the latter is an
isomorphism of $\bR[q]$-modules, but {\it not} of algebras;
following \cite{FGP}, we call it the {\it quantization map}. So
the main goal of our paper is to give a presentation of the
quantization map. For $G=SL(n,\bC)$, the problem has been solved
by Fomin, Gelfand and Postnikov \cite{FGP}. We are going to extend
their result to an arbitrary semisimple  Lie group $G$.

As  in the non-quantum case, we will essentially rely on the
Chevalley formula, this time in its quantum version: the formula
was obtained by D. Peterson in \cite{P} (for more details, see
 section 10 of Fulton and Woodward \cite{FW}).   If $\alpha^{\vee}$ is a
positive  coroot, we consider its {\it height}
$$|\alpha^{\vee}|=m_1+\ldots +m_l,$$ where the positive integers
$m_1, \ldots, m_l$ are given by $\alpha^{\vee}=
m_1\alpha_1^{\vee}+\ldots +m_l\alpha_l^{\vee}$. We also put
$$q^{\alpha^{\vee}}=q_1^{m_1} \ldots q_l^{m_l}.$$ \proclaim{Theorem
3.1} {\rm (Quantum Chevalley Formula; Peterson \cite{P}, Fulton
and Woodward \cite{FW})}
 In $(QH^*(G/B),\circ)$ one has
$$\sigma_{s_i}\circ \sigma_w =\sigma_{s_i}\sigma_w
+ \sum_{l(ws_{\alpha})=
l(w)-2|\alpha^{\vee}|+1}\lambda_i(\alpha^{\vee})q^{\alpha^{\vee}}
\sigma_{ws_{\alpha}}.$$
\endproclaim
The following inequality can be found in Peterson's notes
\cite{P}, as well as in Brenti, Fomin and Postnikov \cite{BFP}.
For the sake of completeness, we will give our own proof of it.
\proclaim{Lemma 3.2}
 For any positive root $\alpha$ we have
$l(s_{\alpha})\leq 2 |\alpha^{\vee}|-1$.
\endproclaim
\demo{Proof} We  prove the lemma by induction on $l(s_{\alpha})$.
If $l(s_{\alpha})=1$, then $\alpha$, as well as $\alpha^{\vee}$,
is simple, so $|\alpha^{\vee}|=1$. Let now $\alpha$ be a positive,
non-simple root. There exists a simple root $\beta$ such that
$\alpha(\beta^{\vee})>0$ (otherwise we would be led to
$\alpha(\alpha^{\vee})\leq 0$). Consequently,
$\beta(\alpha^{\vee})$ is a strictly positive number, too, hence
$$s_{\alpha}(\beta)=\beta-\beta(\alpha^{\vee})\alpha$$ must be a
negative root. Also
$$s_{\beta}s_{\alpha}(\beta)=(\alpha(\beta^{\vee})\beta(\alpha^{\vee})-1)\beta-
\beta(\alpha^{\vee})\alpha$$ is a negative root. By Lemma 3.3,
chapter 1 of \cite{Hi}, we have
$l(s_{\beta}s_{\alpha}s_{\beta})=l(s_{\alpha})-2$. Because
$$s_{\beta}
(\alpha)^{\vee}=s_{\beta}(\alpha^{\vee})=\alpha^{\vee}-\beta(\alpha^{\vee})\beta^{\vee},$$
we have $|s_{\beta}
(\alpha)^{\vee}|=|\alpha^{\vee}|-\beta(\alpha^{\vee})$. By the
induction hypothesis we conclude:
$$l(s_{\alpha})=l(s_{\beta}s_{\alpha}s_{\beta})+2\leq
2|s_{\beta}(\alpha)^{\vee}|-1+2=2|\alpha^{\vee}|-1+2(1-\beta(\alpha^{\vee}))\leq
2|\alpha^{\vee}|-1.$$  \qed
\enddemo
Denote by $\tilde{\Phi}^+$ the set  of all positive  roots
$\alpha$ with the property $l(s_{\alpha})=2|\alpha^{\vee}|-1.$ The
following operators
$$\Lambda_i=\lambda_i +\sum_{\alpha\in \tilde{\Phi}^+}\lambda_i
(\alpha^{\vee})q^{\alpha^{\vee}}\Delta_{s_{\alpha}} \tag 5$$ on
$\bR[q,\lambda]$, $1\leq i \leq l$ have been considered by
Peterson in \cite{P}. His key observation is that we have
$$\Lambda_i[c_w]= \lambda_i[c_w] +
\sum_{l(ws_{\alpha})=
l(w)-2|\alpha^{\vee}|+1}\lambda_i(\alpha^{\vee})q^{\alpha^{\vee}}
[c_{ws_{\alpha}}], \tag 6 $$ the right hand side being, by the
quantum Chevalley formula, just $\lambda_i \circ [c_w]$. In order
to justify (6), we only have to say that if $w\in W$ and $\alpha$
is a positive root with $l(ws_{\alpha})= l(w)-2|\alpha^{\vee}|+1$,
then, by Lemma 3.2, $\alpha$ must be in $\tilde{\Phi}^+$.

 From the
associativity of the quantum product $\circ$ it follows that any
two $\Lambda_i$ and $\Lambda_j$ commute {\it as operators on
$(\bR[\lambda]/I_W)\otimes \bR[q]$}. In fact the following  stronger
result (also stated by Peterson in \cite{P}) holds:
\proclaim{Lemma 3.3} The operators $\Lambda_1, \ldots, \Lambda_l$
on $\bR[q, \lambda]$ commute.
\endproclaim
\demo{Proof} Put $w=s_{\alpha}$ in Lemma 2.2 and obtain:
$$\Delta_{s_{\alpha}}\lambda_i^*=(\lambda_i^*
-\lambda_i(\alpha^{\vee})\alpha^*)\Delta_{s_{\alpha}}
+\sum_{\gamma\in \Phi^+,
l(s_{\alpha}s_{\gamma})=l(s_{\alpha})-1}\lambda_i(\gamma^{\vee})
\Delta_{s_{\alpha}s_{\gamma}}.$$ It follows
$$\align \Lambda_j\Lambda_i=&(\lambda_j\lambda_i)^*
+\sum_{\alpha\in\tilde{\Phi}^+}\lambda_i(\alpha^{\vee})
q^{\alpha^{\vee}}\lambda_j^*\Delta_{s_{\alpha}}\\{}&
+\sum_{\alpha\in\tilde{\Phi}^+}\lambda_j(\alpha^{\vee})
q^{\alpha^{\vee}}\lambda_i^*\Delta_{s_{\alpha}}
-\sum_{\alpha\in\tilde{\Phi}^+}\lambda_j(\alpha^{\vee})
\lambda_i(\alpha^{\vee})q^{\alpha^{\vee}}\alpha^*\Delta_{s_{\alpha}}\\{}&
+\sum_{\alpha\in\tilde{\Phi}^+, \gamma \in \Phi^+,
l(s_{\alpha}s_{\gamma})= l(s_{\alpha})-1}\lambda_j(\alpha^{\vee})
\lambda_i(\gamma^{\vee})q^{\alpha^{\vee}}\Delta_{s_{\alpha}
s_{\gamma}} \\{}& + \sum_{\alpha,\beta \in \tilde{\Phi}^+,
l(s_{\alpha}s_{\beta})=l(s_{\alpha})+l(s_{\beta})}
\lambda_j(\alpha^{\vee})
\lambda_i(\beta^{\vee})q^{\alpha^{\vee}+\beta^{\vee}}
\Delta_{s_{\alpha}s_{\beta}}.
\endalign $$
Denote by $\Sigma_{ij}$ the sum of the last two sums:
the rest is obviously invariant by interchanging $i\leftrightarrow j$.

Let us return to the Bernstein-Gelfand-Gelfand construction
 described in the first section:
Fix $c_{w_0}\in\bR[\lambda]$ such that $[c_{w_0}]=\sigma_{w_0}$
and then set $c_w=\Delta_{w^{-1}w_0}c_{w_0}$, $w\in W$; their
classes modulo $I_W$ are a basis of
 $\bR[ \lambda]/I_W$.
As we said earlier,   from the associativity of the quantum
product we deduce that $\Lambda_j\Lambda_i[c_w]$ is symmetric in
$i$ and $j$, for any $w\in W$. In particular,
$\Sigma_{ij}[c_{w_0}]$ is symmetric in $i$ and $j$. Because
$l(w_0v)=l(w_0)-l(v)$ for any $v\in W$, we have
$$\align \Sigma_{ij}[c_{w_0}] =&\sum_{\alpha\in\tilde{\Phi}^+,
l(s_{\alpha}s_{\gamma})= l(s_{\alpha})-1}\lambda_j(\alpha^{\vee})
\lambda_i(\gamma^{\vee})q^{\alpha^{\vee}}[c_{w_0s_{\gamma}
s_{\alpha}}]\\{}&+ \sum_{\alpha,\beta \in \tilde{\Phi}^+,
l(s_{\alpha}s_{\beta})=l(s_{\alpha})+l(s_{\beta})}
\lambda_j(\alpha^{\vee})
\lambda_i(\beta^{\vee})q^{\alpha^{\vee}+\beta^{\vee}}
[c_{w_0s_{\beta}s_{\alpha}}].
\endalign $$
The latter reproduces exactly the expression of $\Sigma_{ij}$
itself: $\{[c_w]:w \in W\}$ (actually $\{[c_{w_0w^{-1}}]:w\in
W\}$)
 are linearly independent, exactly
like the operators $\{\Delta_w:w\in W\}$. So $\Sigma_{ij}$ is
symmetric in $i$ and $j$ and the lemma is proved.
\qed
\enddemo

The next result is a generalization of Lemma 5.3 of \cite{FGP}.

\proclaim{Lemma 3.4} The map $\psi: \bR[q, \lambda] \to \bR[q,
\lambda]$ given by
$$f\mapsto f(\Lambda_1,\ldots , \Lambda_l)(1)$$
is an $\bR[q]$-linear isomorphism. If $f\in \bR[q,\lambda]$ has
degree $d$ with respect to $\lambda_1,\ldots ,\lambda_l$, then we
can express $\psi^{-1}(f)$ as follows
$$\align \psi^{-1}(f)=&\frac{I-(I-\psi)^{d}}{\psi}(f)\\=&{d\choose 1} f-{{d}\choose 2}\psi(f)
+\ldots + (-1)^{d -2}{{d}\choose{d-1}}\psi^{d
-2}(f)+(-1)^{d-1}\psi^{d -1}(f), \tag 7 \endalign
$$ where ${d\choose 1} ,\ldots , {d \choose {d-1}}$ are the
binomial coefficients.
\endproclaim
\demo{Proof} The degrees of elements of $\bR[q,\lambda]$ we are
going to refer to here are  taken {\it only} with respect to
$\lambda_1,\ldots, \lambda_l$. First, $\psi$ is injective, because
if $g\in \bR[q,\lambda]$ has the property that
$g(\Lambda_1,\ldots,\Lambda_l)(1)=0$, then obviously $g$ must be
$0$. In order to prove both  surjectivity and the formula for
$\psi^{-1}$, we notice that the operator $I-\psi$ lowers the
degree of a polynomial  by at least one, so if $f$ is a
polynomial of degree $d$, then $(I-\psi)^d(f)=0$.\qed
\enddemo

The next result is a direct consequence of the quantum Chevalley
formula.

\proclaim{Proposition 3.5} For any of the generators $R_1,\ldots,
R_l$ of the ideal $I_W^q$, $\psi(R_i)$ is\footnote{In view of
Theorem 5.5 of \cite{FGP}, we could actually expect to have
 $\psi(R_i)=u_i$.} an  $\bR[q]$-linear
combination  of  elements  of $I_W$, the free term with respect to
$q_1, \ldots, q_l$ being $u_i$.
   Hence $\psi(I_W^q)=I_W\otimes \bR[q]$ and $\psi$ gives rise
to a bijection $$\psi : \bR[q,\lambda]/I_W^q\to \bR[q, \lambda]/(I_W\otimes
\bR[q]).$$
\endproclaim
\demo{Proof} We just have to use the fact that
$$\lambda_{i_1}\circ\ldots \circ \lambda_{i_k}=
\Lambda_{i_1} \ldots  \Lambda_{i_k} (1)  \text{ mod} \  I_W\otimes
\bR[q]$$ so that
$$\align \psi(R_i)  \text{ mod} \   I_W\otimes \bR[q]&
=R_i(q_1,\ldots, q_l,\Lambda_1,\ldots, \Lambda_l)(1)  \text{ mod}
\
 I_W\otimes \bR[q]\\{} &= R_i (q_1,\ldots, q_l,\lambda_1\circ ,\ldots,
\lambda_l\circ)\\{}&=0. \endalign$$\qed
\enddemo

Our  polynomial representatives of  Schubert classes in
$QH^*(G/B)$ are described by the following theorem, which is the
central result of the paper. The proof is governed by the same
ideas that have been used in the non-quantum case (see section 2).
\proclaim{Theorem 3.6} The quantization map $\bR[q, \lambda]/(I_W
\otimes \bR[q])\to \bR[q,\lambda]/I_W^q$ is just $\psi^{-1}$. More
precisely, if $w\in W$ has length $l(w)=l$, then the class of
$c_w$ in $\bR[q, \lambda]/(I_W\otimes \bR[q])$ is mapped to the
class of
$$ \frac{I-(I-\psi)^{l}}{\psi}(c_w)
={l\choose 1} c_w-{{l}\choose 2}\psi(c_w) +\ldots + (-1)^{l
-2}{{l}\choose{l-1}}\psi^{l -2}(c_w)+(-1)^{l-1}\psi^{l -1}(c_w)
$$ in $\bR[q, \lambda]/I_W^q$, where $\psi$ has been defined in Lemma
3.4.
\endproclaim
\demo{Proof} For any polynomial $f\in\bR[q,\lambda]$, we denote by
$[f]$, $[f]_q$ its classes modulo $I_W\otimes \bR[q]$,
respectively modulo $I_W^q$. By the definition of $\psi$, the
polynomial $\hat{c}_w:=\psi^{-1}(c_w)$ is determined by
$$\hat{c}_w(\Lambda_1,\ldots,\Lambda_l)(1)=c_w.$$  We take into
account  (6), where $\Lambda_i[c_w]$ is the same as
$$[\Lambda_i(c_w)]= [\Lambda_i(\hat{c}_w(\Lambda_1, \ldots,
\Lambda_l)(1))] =\psi ([\lambda_i\hat{c}_w]_q).$$ Because
$[c_v]=\psi([\hat{c}_v]_q)$ for any $v\in W$ and the map $\psi$ is
bijective, it follows that in $\bR[q, \lambda]/I_W^q$ we have
$$\lambda_i [\hat{c}_w]_q= \sum_{l(ws_{\alpha})=
l(w)+1}\lambda_i(\alpha^{\vee}) [\hat{c}_{ws_{\alpha}}]_q +
\sum_{l(ws_{\alpha})=
l(w)-2|\alpha^{\vee}|+1}\lambda_i(\alpha^{\vee})q^{\alpha^{\vee}}
[\hat{c}_{ws_{\alpha}}]_q. \tag 8$$

As $\bR[q]$-algebras, both $QH^*(G/B)$ and $\bR[q, \lambda]/I_W^q$
are generated by their degree 2 elements; this is why their
structure is uniquely determined by the bases $\{\sigma_w:w\in
W\}$, respectively $\{[\hat{c}_w]:w\in W\}$ and the matrices of
multiplication
 by $\sigma_{s_i}$, respectively $\lambda_i$, $1 \leq i \leq l$.
Since $\hat{c}_{s_i}=\lambda_i$, $1\leq i \leq l$, it follows from
Theorem 3.1 and   relation (8) that the map $$QH^*(G/B)\to \bR[q,
\lambda]/I_W^q \text{ \ given \ by \ } \sigma_w\mapsto \hat{c}_w,
w\in W$$ is an isomorphism of algebras and the proof is
finished.\qed
\enddemo

\subheading{Example} We will illustrate our main result by giving
concrete solutions to the quantum Giambelli problem for $G/B$,
where $G$ is simple of type $B_2$. This is the first interesting
case, different from $A_n$ and for which $\tilde{\Phi}^+\neq
\Phi^+$. We will use the following presentation of the root
system: if $x_1, x_2$ are an orthogonal coordinate system of the
plane and $e_1, e_2$ the unit direction vectors of the coordinate
axes, then

\item - the simple roots are $\alpha_1:=x_1$ and
$\alpha_2:=x_2-x_1$.

\item - the positive roots are $\alpha_1, \alpha_2, \alpha_3:=
\alpha_1+\alpha_2=x_2$ and $\alpha_4:=2\alpha_1+\alpha_2=x_1+x_2$.

\item - the  positive coroots are $\alpha_1^{\vee}=2e_1$,
$\alpha_2^{\vee}=e_2-e_1$,
$\alpha_3^{\vee}=2e_2=\alpha_1^{\vee}+2\alpha_2^{\vee}$ and
$\alpha_4^{\vee}=e_1+e_2=\alpha_1^{\vee}+\alpha_2^{\vee}$.

\item - the fundamental weights $\lambda_1, \lambda_2$ are
determined by
$$\align {}& x_1=2\lambda_1-\lambda_2 \\
{}& x_2=\lambda_2
\endalign$$

\item - the simple reflections are $s_1: (x_1,x_2)\mapsto (-x_1,
x_2)$ and $s_2: (x_1, x_2)\mapsto (x_2, x_1)$. The generators of
$I_W$ are obviously $x_1^2+x_2^2$ and $x_1^2x_2^2$.

\item - following \cite{FK}, we can obtain  polynomial
representatives of Schubert classes in $\bR[x_1,x_2]/(x_1^2+x_2^2,
x_1^2x_2^2)$ as indicated in the following table:

\newpage

\bigpagebreak

\settabs 4\columns

\+  & $w$ & $c_w$ &\cr

\+ & -------------------- &
-------------------------------- &\cr

\+ & $w_0=s_1s_2s_1s_2$ & $(x_1-x_2)^3(x_1+x_2)/16$ &\cr

\+ & $s_2s_1s_2$ & $-x_2(x_1-x_2)(x_1+x_2)/4$ &\cr

\+ & $s_1s_2s_1$ & $-(x_1-x_2)^2(x_1+x_2)/8$ &\cr

\+ & $s_2s_1$ & $(x_1+x_2)^2/4$ &\cr

\+ & $s_1s_2$ & $-(x_1-x_2)(x_1+x_2)/4$ &\cr

\+ & $s_2$ & $x_2$ &\cr

\+ & $s_1$ & $(x_1+x_2)/2$ &\cr

\medpagebreak

Note that we have started the B-G-G algorithm with $c_{w_0}$
which differs from
$\alpha_1\alpha_2\alpha_3\alpha_4/8$  by a multiple of $x_1^2+x_2^2$.

Theorem 2.6 will allow us to describe the quantization map without
knowing anything about  the ideal $I_W^q$ of quantum relations.
But for the sake of completeness we will also obtain the two
generators of $I_W^q$, by using the theorem of Kim as presented in
our paper \cite{M}. We have to consider  the Hamiltonian system
which consists of the standard 4-dimensional symplectic manifold
$(\bR^4, dr_1\wedge ds_1+dr_2\wedge ds_2)$ with the Hamiltonian
function
$$E(r,s)=\sum_{i,j=1}^2\langle
\alpha_i^{\vee},\alpha_j^{\vee}\rangle
r_ir_j+\sum_{i=1}^2e^{-2s_i}=(2r_1-r_2)^2+r_2^2+e^{-2s_1}+e^{-2s_2}.$$
The first integrals of motion of the system are $E$ and --- by
inspection
--- the function
$$F(r,s)=
(2r_1-r_2)^2r_2^2+r_2^2e^{-2s_1}-(2r_1-r_2)r_2e^{-2s_2}+2e^{-2s_1}e^{-2s_2}+
\frac{1}{4}(e^{-2s_2})^2.$$ By the main result of \cite{M}, the
quantum relations are obtained from $E$, respectively $F$, by the
formal replacements:
$$\align {}& 2r_1-r_2 \mapsto x_1, r_2\mapsto x_2  \\
{}& e^{-2s_1}\mapsto -\langle \alpha_1^{\vee},
\alpha_1^{\vee}\rangle q_1=-4q_1, e^{-2s_2}\mapsto -\langle
\alpha_2^{\vee}, \alpha_2^{\vee}\rangle q_2=-2q_2. \endalign $$ In
conclusion, $I_W^q$ is the ideal of $\bR[q_1,q_2,x_1,x_2]$
generated by
$$ x_1^2 +x_2^2 -4q_1-2q_2=0 \ \text{and } \
 x_1^2 x_2^2-4q_1 x_2^2+2q_2 x_1 x_2+16 q_1q_2+q_2^2.
 $$

 Now, we will
determine explicitly the image of each Schubert class $\sigma_w$,
$w\in W$ via the isomorphism
$$QH^*(G/B)\simeq \bR[q_1,q_2,x_1, x_2]/I_W^q.$$
  The place of the
operators $\Lambda_1, \Lambda_2$ is taken by $X_1, X_2$ where
$$X_i=x_i+x_i(\alpha_1^{\vee})q_1\Delta_{s_1}+x_i(\alpha_2^{\vee})q_2\Delta_{s_2}+
x_i(\alpha_4^{\vee})q_1q_2\Delta_{s_1}\Delta_{s_2}\Delta_{s_1},\qquad
i=1,2.
$$ More precisely, we have
$$X_1=x_1+2q_1\Delta_{s_1}-q_2\Delta_{s_2}+q_1q_2\Delta_{s_1}\Delta_{s_2}\Delta_{s_1}$$
and
$$X_2=x_2+q_2\Delta_{s_2}+q_1q_2\Delta_{s_1}\Delta_{s_2}\Delta_{s_1}.$$

Rather than using the formula for $\psi^{-1}$ given by (7), it
seems more convenient to determine $\hat{c}_w:=\psi^{-1}(c_w)\in
\bR[q_1,q_2,x_1,x_2]$ by the definition of $\psi$, i.e. from the
condition
$$\hat{c}_w(X_1, X_2)(1)=c_w(x_1,x_2).$$
We will explain the details just for the case $w=w_0$, which is
the most illustrative one. The polynomial we are looking for has
the form
$\hat{c}_{w_0}=c_{w_0}+q_1a_1+q_2a_2+b_1q_1^2+b_2q_2^2+b_3q_1q_2,$
where $a_1, a_2$ are homogeneous polynomials of degree 2 in $x_1,
x_2$ and $b_1,b_2, b_3$ are constant. The condition that
determines $a_1,a_2, b_1, b_2,b_3$ is
$$\align
c_{w_0}(X_1,X_2)(1)+q_1a_1(X_1,X_2)(1)+q_2a_2(X_1,X_2)(1)+b_1q_1^2+b_2q_2^2+b_3q_1q_2&{}\\
{}=c_{w_0}(x_1,x_2).\tag 9 \endalign$$

 The first step is  to compute $c_{w_0}(X_1,X_2)(1)$ and determine $a_1$ and $a_2$.
  Using 
$$\Delta_{s_1}(f)=\frac{f(x_1,x_2)-f(-x_1,x_2)}{x_1} \ \  \text{
and } \  \
\Delta_{s_2}(f)=\frac{f(x_1,x_2)-f(x_2,x_1)}{x_2-x_1},$$ $f\in
\bR[x_1,x_2]$ we obtain
$$c_{w_0}(X_1,X_2)(1)=c_{w_0}(x_1,x_2)+\frac{1}{8}q_1(3x_1^2-4x_1x_2+x_2^2)+
\frac{1}{4}q_2(x_1^2-x_2^2)+q_1^2+q_1q_2.$$ Since the coefficients
of $q_1$, respectively $q_2$ in the left hand side of (9) must
vanish, we deduce:
$$a_1=-\frac{1}{8}(3x_1^2-4x_1x_2+x_2^2), \ \ \
a_2=-\frac{1}{4}(x_1^2-x_2^2).$$

The second step is to compute $a_1(X_1,X_2)(1)$ and
$a_2(X_1,X_2)(1)$ and determine $b_1, b_2$ and $b_3$.  We take
into account that
$$X_1-X_2=x_1-x_2+2q_1\Delta_{s_1}-2q_2\Delta_{s_2}$$  and
find
$$\align {}&a_1(X_1,X_2)(1)=-\frac{1}{8}(X_1-X_2)(3x_1-x_2)=a_1(x_1,x_2)-\frac{3}{2}q_1-q_2\\
{}&a_2(X_1,X_2)(1)=-\frac{1}{4}(X_1-X_2)(x_1+x_2)=a_2(x_1,x_2)-q_1.\endalign$$
Coming back to (9), we deduce
$$b_1=\frac{1}{2}, \  b_2=0, \  b_3=1,$$
hence
$$\hat{c}_{w_0}=c_{w_0}-\frac{1}{8}q_1(3x_1^2-4x_1x_2+x_2^2)-\frac{1}{4}q_2(x_1^2-x_2^2)
+\frac{1}{2}q_1^2+q_1q_2.$$

The other $\hat{c}_w$, $w\in W$, can be obtained by similar
computations. They are described in the following table:

\bigpagebreak

\settabs 4\columns

\+  & $w$ & $\hat{c}_w-c_w$ &\cr

\+ & -------- & -------------------------------- &\cr

\+ & $s_2s_1s_2$ & $q_1x_2$ &\cr

\+ & $s_1s_2s_1$ &
$\frac{1}{2}(x_1-x_2)q_1+\frac{1}{2}(x_1+x_2)q_2$ &\cr

\+ & $s_2s_1$ & $-q_1$ &\cr

\+ & $s_1s_2$ & $q_1$ &\cr

\+ & $s_2$ & $0$ &\cr

\+ & $s_1$ & $0$ &\cr

\medpagebreak

\subheading{Acknowledgements}  I would like to thank Martin Guest
and Takashi Otofuji for the extensive exchange of ideas which led
me to Theorem 3.6. I am also grateful  to McKenzie Wang and Chris
Woodward for suggesting  improvements to previous versions of the
paper. Finally, I would like to thank the referee for several
helpful suggestions.

\bigskip

\eightpoint \it

\noindent  Department of Mathematics and Statistics
\newline
McMaster University
\newline
1280 Main Street West
\newline
Hamilton, Ontario, Canada
\newline
L8S 4K1

\bigskip

\noindent \rm{and} 

\it

\bigskip

\noindent  Department of Mathematics
\newline
 University of Toronto
\newline
100 St. George Street 
\newline
Toronto, Ontario, Canada
\newline
M5S 3G3

\noindent
amare\@math.toronto.edu

\enddocument